\documentclass{PJM}

\authorone{Thando Nkomozake}
\addressone{Research and Development Department, Vision Tech Strategies,
Bellville, Cape Town, 7530}
\countryone{South Africa}
\emailone{thandolwethu3@gmail.com }
\title{Curve Closest to Sphere}
\headlinetitle{Curve Closest To Sphere}
\begin{document}
\vskip0.2in
\vskip0.1in
\normalsize

 {\bf Abstract} We propose a solution to the tenth of Professor Clark Kimberling's unsolved problems found on https://faculty.evansville.edu/ck6/integer/unsolved.html. We are required to find the parametric equations of a simple and closed curve $C$ on the unit sphere $S$ with arc-length $4 \pi$, that minimizes the mean arc-distance from $S$ to $C$. We give explicit definitions of the mean arc-distance from $C$ to $S$, $M$ and the mean arc-distance from $S$ to $C$, $\tilde{M}$. We show that these two quantities are not the same. We show that for all closed and simple curves $C$ of arc-length $4 \pi$ on $S$, $M$ is constant and is equal to $2 \pi^{2}$. Therefore all such curves minimize $M$. We show that in contrast, $\tilde{M}$ varies for different closed and simple curves $C$ of arc-length $4 \pi$ on $S$. We find such a curve that minimizes $\tilde{M}$.

\maketitle

\section{Mean Arc-Distance From $C$ to $S$}       

\noindent Let $\Vec{x} = (x, y, z)$ be an arbitrary point on $S$ and let $ \Vec{r} (t) $ be a parameterization of a simple and closed curve $C$ that lies on $S$. Then the mean arc-distance from $ \Vec{r} (t) $ to $S$ is given by 
\begin{equation}
    \Bar{D} = \frac{1}{4 \pi} \int_{S}  \textbf{dist} (\Vec{r} (t), \Vec{x}) \ dS, 
\end{equation}
where $ \textbf{dist} (\Vec{r} (t), \Vec{x}) = \arccos\big(  \sin(\theta(t)) \cos(\phi(t)) \ x + \sin(\theta(t)) \sin(\phi(t)) \ y + \cos(\theta(t)) \ z
  \big) $, which of course gives the shortest distance between $\Vec{x}$ and $\Vec{r} (t)$ on $S$. So we are required to find $C$-as defined above and with parameterization $\Vec{r} (t)$-that minimizes $\Bar{D}$. Since $\Bar{D}$ is a surface integral over $S$ we can calculate the surface integral first before we minimize. \\
  Let $ A = \sin(\theta(t)) \cos(\phi(t))$, $ B = \sin(\theta(t)) \sin(\phi(t))$ and $ E = \cos(\theta(t))$. Then we get 
\begin{equation}
    \Bar{D} = \frac{1}{4 \pi} \int_{S}   \arccos\big(  A x + B  y + E z
  \big) \ dS.
\end{equation}We then parameterize the integral $\Bar{D}$ by letting $x =\sin(\theta) \cos(\phi) $, $y = \sin(\theta) \sin(\phi) $ and $z = \cos(\theta)$, where $0 \leq \theta \leq \pi $ and $0 \leq \phi \leq 2 \pi$. This implies \begin{equation}
    \Bar{D} = \frac{1}{4 \pi} \int_{0}^{2 \pi} \int_{0}^{\pi} 
\arccos\big(  A \sin(\theta) \cos(\phi) + B \sin(\theta) \sin(\phi) +  E\cos(\theta)  
  \big) \ |\Vec{r}_{\theta} \times \Vec{r}_{\phi} | \ d \theta \ d \phi .
\end{equation}\\ Now we let $D = A \cos(\phi) + B \sin(\phi)$ which implies \begin{equation}
    \arccos\big(  A \sin(\theta) \cos(\phi) + B \sin(\theta) \sin(\phi) + E \cos(\theta)\big) = \arccos\big(D \sin(\theta) + E \cos(\theta)
  \big).
\end{equation} By the identity relating arccos and arcsin [1] we have \begin{equation}
    \sin(\theta) \arccos\big(D \sin(\theta) + E \cos(\theta)
  \big) = \frac{\pi}{2} \sin(\theta) - \sin(\theta) \arcsin\big(D \sin(\theta) + E \cos(\theta)
  \big).
\end{equation}Using the power series of arcsin [1] we get 
\begin{equation}
\begin{split}
    \sin(\theta) \arcsin\big(D \sin(\theta) + E \cos(\theta)
  \big) & = \sum_{n \geq 0} \sum_{k = 0}^{2 n + 1} 
  \frac{ 1 }{ (2n + 1) 4^{n}} \ {2n \choose n} {{2n + 1} \choose k } \\
  & \times  D^{k} \ E^{2n + 1-k} \ \sin^{k + 1}(\theta) \ \cos^{2n +1-k}(\theta).
 \end{split}
 \end{equation}

\noindent We then have \begin{equation}
    \begin{split}
        \int_{0}^{2 \pi} \int_{0}^{\pi} \ \sin(\theta) \arccos\big(D \sin(\theta) + E \cos(\theta)
  \big) \ d \theta \ d \phi &= \\
  \int_{0}^{2 \pi} \int_{0}^{\pi} \ \frac{\pi}{2} \sin(\theta) \ d \theta \ d \phi - \int_{0}^{2 \pi} \int_{0}^{\pi} \ \sin(\theta) \arcsin\big(D \sin(\theta) + E \cos(\theta)\big) \ d \theta \ d \phi,
    \end{split}
\end{equation} \\

\noindent where 
\begin{equation}
    \begin{split}
        \int_{0}^{2 \pi} \int_{0}^{\pi} \ \sin(\theta) \arcsin\big(D \sin(\theta) + E \cos(\theta)\big) \ d \theta \ d \phi &= \\
        \sum_{n \geq 0} \sum_{k = 0}^{2 n + 1} 
  \frac{ 1 }{ (2n + 1) 4^{n}} \ {2n \choose n} {{2n + 1} \choose k } \\
  \times \int_{0}^{\pi} \ \sin^{k+1}(\theta) \ \cos^{2n + 1 - k}(\theta) \ \bigg(  \int_{0}^{2 \pi}  D^{k} \ E^{2n + 1 - k}  \ d\phi    \bigg) \ d\theta.
    \end{split}
\end{equation}\\
 Using the binomial theorem [2] we get \begin{equation}
     \int_{0}^{2 \pi}  D^{k} \ E^{2n + 1 - k}  \ d\phi = 0,
 \end{equation}
which then implies \begin{equation}
    \int_{0}^{2 \pi} \int_{0}^{\pi} \ \sin(\theta) \arcsin\big(D \sin(\theta) + E \cos(\theta)\big) \ d \theta \ d \phi = 0.
\end{equation}\\
We also have \begin{equation}
    \int_{0}^{2 \pi} \int_{0}^{\pi} \ \frac{\pi}{2} \sin(\theta) \ d \theta \ d \phi = 2 \pi^{2},
\end{equation}\\
which then implies  \begin{equation}
    \int_{0}^{2 \pi} \int_{0}^{\pi} \ \sin(\theta) \arccos\big(D \sin(\theta) + E \cos(\theta)
  \big) \ d \theta \ d \phi = 2 \pi^{2}.
\end{equation}\\
This then implies $ \Bar{D} = \frac{1}{4 \pi} \big( 2 \pi^{2} \big) = \frac{\pi}{2}$.\\

\noindent So this shows that the mean arc-distance from an arbitrary point $\Vec{r}(t)$ on a curve $C$ that lies on $S$ to $S$ is the constant $\frac{\pi}{2}$. We can then define the mean arc-distance from a curve $C$ as a whole to $S$ as \begin{equation}
    \begin{split}
        M &= \int_{C} \ \frac{\pi}{2} \ ds\\
        &= \int_{t_{i}}^{t_{f}} \ \frac{\pi}{2} \ |\Vec{r} \ ^\prime(t)| \ dt \\
        &= \frac{\pi}{2} \int_{t_{i}}^{t_{f}}  \ |\Vec{r} \ ^\prime(t)| \ dt.
    \end{split}
\end{equation}
Therefore minimizing $M$ would be minimizing the arc-length of $C$. However if we state the condition that the curve $C$ must have arc-length $4 \pi$, that means that all curves $C$ with this same arc-length will give the same $M = \frac{\pi}{2} \big(  4 \pi \big) = 2 \pi^{2} $. Meaning that they all minimize $M$ and the minimum they give is $2 \pi^{2}$. 

\section{Mean Minimum Arc-Distance}

\noindent Alternatively you may be compelled to approach this problem in the following way: suppose that we start with an arbitrary closed and simple curve $C$ with an arc-length $4 \pi$ on $S$ and a parameterization $ \Vec{r} (t) $. Let $ \Vec{x}_{1}$ be an arbitrary point on $S$. Then the arc-distance from $ \Vec{x}_{1}$ to a point $ \Vec{r} (t) $ on $C$ is given by $\arccos\big( \Vec{r} (t) \cdot \Vec{x}_{1}   \big)$. In order to minimize this arc-distance we would have to find the $t$ value, say $t_{1}$, for which this arc-distance is the minimum. We would then have that the minimum arc-distance from $ \Vec{x}_{1}$ to $C$ is given by $\arccos\big( \Vec{r} (t_{1}) \cdot \Vec{x}_{1}   \big)$. We would do this for other points on $S$ say $\Vec{x}_{1}, \Vec{x}_{2},..., \Vec{x}_{n}$. Then to calculate the approximate mean minimum arc-distance from $S$ to $C$ we sum as follows \begin{equation}
     \frac{1}{4 \pi}\bigg[ \arccos\big( \Vec{r} (t_{1}) \cdot \Vec{x}_{1}   \big) + \arccos\big( \Vec{r} (t_{2}) \cdot \Vec{x}_{2}   \big) + ... + \arccos\big( \Vec{r} (t_{n}) \cdot \Vec{x}_{n}   \big)           \bigg].
\end{equation} Then we would have to use this sum (2.1) to find the curve $C$ that minimizes it. Note that the sum (2.1) above is the mean minimum arc-distance from $S$ to $C$, which is different from the mean arc-distance from $S$ to $C$. This approach is not feasible because in order to find ${t_{1}}, {t_{2}},..., {t_{n}}$ we need $ \Vec{r} (t) $ which is also unknown. It also deals with the mean minimum arc-distance, which is different from the mean arc-distance that the problem statement is referring to.\\

\section{Mean Arc-Distance From $S$ to $C$}

\noindent The appropriate alternative approach is as follows: we start with an arbitrary closed and simple curve $C$ of arc-length $4 \pi$ on $S$. Let the parameterization of $C$ be given by\\
$\Vec{r} (t) = \big( \sin(\theta(t)) \cos(\phi(t)), \sin(\theta(t)) \sin(\phi(t)), \cos(\theta(t))              \big)$, where $t \in [t_{i}, t_{f}]$. Let $P$ be an arbitrary point on $S$ with coordinates 
$\big( \sin( \theta_{0} ) \cos( \phi_{0} ), \sin( \theta_{0} ) \sin( \phi_{0} ), \cos( \theta_{0} )        \big)$. Let $Q$ be an arbitrary point on $C$ with coordinates $\Vec{r} (t)$. Then the arc-distance from $P$ to $Q$ is given by\\
$  \arccos\big( \cos(\theta_{0}) \cos(\theta(t))  + \sin(\theta_{0}) \sin(\theta(t)) \cos(\phi_{0}-\phi(t)) 
  \big) $. We now calculate the mean arc-distance from $P$ to $C$ by integrating the arc-distance from $P$ to $Q$ over $[t_{i}, t_{f}]$. This gives \begin{equation}
      \tilde{S} = \frac{1}{t_{f} - t_{i}} \int_{t_{i}}^{t_{f}} \ \arccos\big(\cos(\theta_{0}) \cos(\theta(t))  + \sin(\theta_{0}) \sin(\theta(t)) \cos(\phi_{0}-\phi(t)) 
  \big) \ dt.
  \end{equation} We now find the curve $C$ that minimizes $\tilde{S}$. We let \\ 
  $L = \arccos\big(\cos(\theta_{0}) \cos(\theta(t))  + \sin(\theta_{0}) \sin(\theta(t)) \cos(\phi_{0}-\phi(t)) 
  \big) $. Then we write the corresponding Euler-Lagrange equations for this problem which are given by \begin{align}
    \frac{\partial L }{\partial \theta} - \frac{d}{dt}\left(\frac{\partial L}{\partial \dot{\theta}}\right) &= 0 \label{eq:el1} ,\\
    \frac{\partial L}{\partial \phi} - \frac{d}{dt}\left(\frac{\partial L}{\partial \dot{\phi}}\right) &= 0. \label{eq:el2}
    \end{align} Because $L$ doesn't contain derivatives these equations reduce to 
    \begin{align}
    \frac{\partial L }{\partial \theta}  &= 0 \label{eq:el1} ,\\
    \frac{\partial L}{\partial \phi}  &= 0. \label{eq:el2}
    \end{align} These equations then give 
    \begin{align}
     \sin(\theta_{0}) \cos(\theta(t)) \cos(\phi_{0}-\phi(t))  &= \cos(\theta_{0}) \sin(\theta(t)) \label{eq:el1} ,\\
    \sin(\theta_{0}) \sin(\theta(t)) \sin(\phi_{0}-\phi(t))  &= 0. \label{eq:el2}
    \end{align} 
    Solving these equations gives the following solution $\theta(t) = m \pi$ and $\phi(t) = \phi_{0} - \big( \frac{\pi}{2} + k \pi \big)$, where $m$ and $k$ $\in \mathbb{Z} $. This solution is discrete because our Euler-Lagrange equations reduced to a system of trigonometric equations. We therefore use this discrete solution to find a closed curve $C$ with arc-length $4 \pi$ that minimizes $\tilde{S}$. Such a curve is the great circle $\tilde{C}:\Vec{r} (t) = \big( \sin(2 \pi t), 0, \cos(2 \pi t)  \big) $ for $t \in [0, 2]$. We substitute this curve into $\tilde{S}$ and get the following integral
    \begin{align}
    \tilde{S} = \frac{1}{2} \int_{0}^{2} \ \arccos\big( \sin(\theta_{0}) \cos(\phi_{0}) \sin(2 \pi t) + \cos(\theta_{0}) \cos(2 \pi t)  \big) \ dt  &= \frac{\pi}{2}.  \label{eq:el1} 
    \end{align} So this means that the mean arc-distance from $P$ to the great circle $\tilde{C}:\Vec{r} (t) = \big( \sin(2 \pi t), 0, \cos(2 \pi t)  \big) $ is $\frac{\pi}{2}$. It's important to note that for any great circle on the unit sphere, the integral (3.8) will give the same result. It's well known that on the unit sphere, the closest closed curve of arc-length $4 \pi$, to an arbitrary point $P$ is a great circle and it gives $\tilde{S} = \frac{\pi}{2}$. Hence $\frac{\pi}{2} $ is the minimum value for $\tilde{S}$ that any closed curve of arc-length $4 \pi$ can give on the unit sphere. Which means that for any closed curve of arc-length $4 \pi$ on the unit sphere we have that $\tilde{S} \geq \frac{\pi}{2}$.\\
    One of the conditions that the desired curve must satisfy is that it must also be simple. In this case the great circle on the unit sphere is not simple because you have to traverse it twice in order for its arc-length to be $4 \pi$. Therefore using numerical methods we found that the curve that satisfies all the required conditions and gives $\tilde{S} \approx \frac{\pi}{2}$, is the tennis ball seam curve shown in the figure 1 below with the following parametric equations (where $A = 0.7037$)
     \begin{equation}
    \begin{split}
       x(t) &= \sin\bigg( \frac{\pi}{2} - \bigg(\frac{\pi}{2} - A \bigg) \cos(t) \bigg) \cos\bigg( \frac{t}{2}  + A \sin(2t) \bigg) \\
       y(t) &= \sin\bigg( \frac{\pi}{2} - \bigg(\frac{\pi}{2} - A \bigg) \cos(t) \bigg) \sin\bigg( \frac{t}{2}  + A \sin(2t) \bigg) \\
       z(t) &= \cos\bigg( \frac{\pi}{2} - \bigg(\frac{\pi}{2} - A \bigg) \cos(t) \bigg).
    \end{split}
\end{equation}\\
This proves that this tennis ball seam curve (3.9) is the desired curve. \\

\begin{figure}[ht!]
\includegraphics[width=1.0\textwidth]{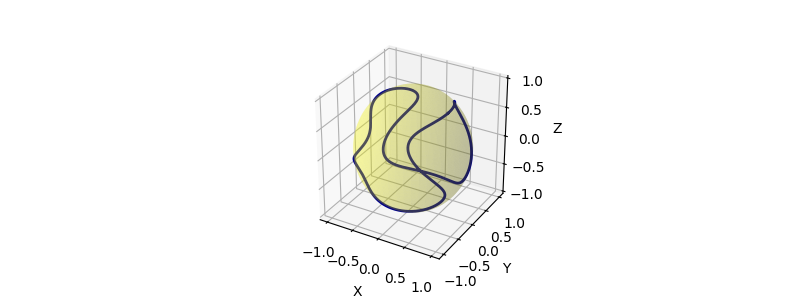} 
\caption{}
  \centering
  
 \end{figure}

\noindent Therefore in order to calculate the mean arc-distance from $S$ to the curve (3.9), we integrate this $\frac{\pi}{2}$ over $S$. This gives 
    \begin{align}
    \tilde{M} = \iint_{S} \ \frac{\pi}{2}  \ dS &= 2 \pi^{2}.  \label{eq:el1} 
    \end{align} Therefore we have shown that the closed and simple curve $C$ with arc-length $4 \pi$, that minimizes the mean arc-distance from $S$ to $C$ is the tennis ball seam curve (3.9) and the minimum that it gives is $\tilde{M} = 2 \pi^{2} $.\\
    
    \noindent We must note that the mean arc-distance from $C$ to $S$, $M$ is a different quantity from the mean arc-distance from $S$ to $C$, $\tilde{M}$. However we can see above that for the curve (3.9), $M = \tilde{M} = 2 \pi^{2} $. In general $M$ and $\tilde{M}$ don't always have the same value. To give an example of this consider the following curve shown in the figure 2 below with parametric equations (where $B = 0.1856$)
     \begin{equation}
    \begin{split}
       x(t) &= \sin\bigg( \frac{3 \pi}{4} + B \sin(10t) \bigg) \cos(t) \\
       y(t) &= \sin\bigg( \frac{3 \pi}{4} + B \sin(10t)  \bigg) \sin(t) \\
       z(t) &= \cos\bigg( \frac{3 \pi}{4} + B \sin(10t)      \bigg).
    \end{split}
\end{equation}\\

\begin{figure}[ht!]
\includegraphics[width=1.0\textwidth]{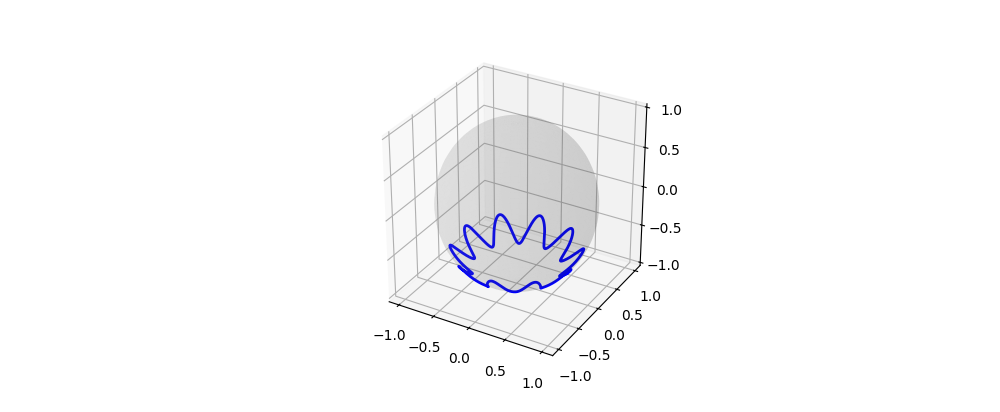} 
\caption{}
  \centering
  
 \end{figure}    

\noindent Let's calculate the mean arc-distance from the arbitrary point $P$ to this curve  (3.11). We get \begin{equation}
    \begin{split}
        \tilde{S} = \frac{1}{2 \pi} \int_{0}^{2 \pi} \ \arccos\bigg( \sin(\theta_{0}) \cos(\phi_{0}) \sin\bigg( \frac{3 \pi}{4} + B \sin(10t) \bigg) \cos(t)                   \\
        + \sin(\theta_{0}) \sin(\phi_{0}) \sin\bigg( \frac{3 \pi}{4} + B \sin(10t)  \bigg) \sin(t)      + \cos(\theta_{0}) \cos\bigg( \frac{3 \pi}{4} + B \sin(10t)      \bigg)        \bigg) \ dt.
    \end{split} 
    \end{equation} The value of this integral in (3.12) varies depending on the point $(\theta_{0}, \phi_{0})$ on $S$. For example if \\
    $\theta_{0} = 0$ and $\phi_{0} = 1$ the integral attains the value $\tilde{S} \approx 2.3562$, which is bigger than $ \frac{\pi}{2} $. This implies that the corresponding mean arc-distance from $S$ to the curve (3.12), $\tilde{M} $ is bigger than $2 \pi^{2}$. So (3.12) is an example of a closed and simple curve $C$ of arc-length $4 \pi$ on $S$, where the mean arc-distance from $S$ to $C$ is not equal to $2 \pi^{2}$. Where in contrast, we've shown that the mean arc-distance from any point on the unit sphere to $S$ is the constant $\Bar{D} = \frac{\pi}{2}  $. Consequently we've shown that for any closed and simple curve $C$ of arc-length $4 \pi$ on $S$, the mean arc-distance from $C$ to $S$ is the constant $M = 2 \pi^{2} $. Which further shows that the mean arc-distance from $C$ to $S$, $M$ is a different quantity from the mean arc-distance from $S$ to $C$, $\tilde{M}$.\\

 \section*{Acknowledgements} 
I would like to thank Professor Clark Kimberling for taking the time to correspond with me. Our conversations were valuable. I would also like to acknowledge and thank the referee from the Palestinian Journal of Mathematics, whose comments were useful and helped me improve the quality of this article. I'm grateful to the Research and Development Department of Vision Tech Strategies for supporting this work.

\end{document}